\documentclass {siamltex} 

\usepackage {amsfonts}
\usepackage[round]{natbib}
\usepackage {verbatim}
\usepackage[T1]{fontenc}

\newcommand {\Vector} {\mbox {\rm vec}}

\newcommand {\JxAvector} {{\mathcal {J}}_x(\Vector (A))}
\newcommand {\JxA} {{\mathcal {J}}_x(A)}

\newcommand {\JxAb} {{\mathcal {J}}_x(A, b)}
\newcommand {\Jxb} {{\mathcal J}_x(b)}
\newcommand {\XxA} {\chi_x(A)}
\newcommand {\XxAb} {\chi_x(A, b)}
\newcommand {\Xxb} {\chi_x(b)}

\renewcommand {\arccos} {\mathop {\mbox {\rm Arccos}} \kern 0.1em}
\newcommand {\bfA} {{\bf A}}
\newcommand {\bfb} {{\bf b}}

\newcommand {\dA} {\delta \kern -0.125em A}
\newcommand {\DA} {\Delta A}

\newcommand {\Db} {\Delta b}

\newcommand {\Dx} {\Delta x}
\newcommand {\reals} {{\mathbb R}}
\newcommand {\bhp} {{\mathcal B}}

\renewcommand {\O} {{\mathcal O}}
\newcommand {\R} {{\mathcal R}}
\newcommand {\Df} [1] {D^{#1} f}
\newcommand {\Dg} [1] {D^{#1} g}
\newcommand {\Oe} {\O (\epsilon^2)}

\newcommand {\Hessian} {{\mathcal H}}
\newcommand {\halfeta} {{\textstyle {\eta \over 2}}}
\newcommand {\param} {\alpha}

\newcommand {\J} {\mathcal J}
\newcommand {\M} {\mathcal M}
\newcommand {\calS} {\mathcal S}
\newcommand {\DM} {\Delta M}
\newcommand {\Dv} {\Delta v}
\newcommand {\plethorasup} [1] {{}^{\mbox {({\it #1}\,)}}}
\newcommand {\plethora} [1] {({\it #1}\/)}
\newcommand {\Sremainder} {{\| \calS \|_2 \over \| x \|_2}}

\newcommand {\fourthpi} {\mbox {\raisebox {0.25ex} {$\pi \over 4$}}}

\newcommand {\onebyone} [1] {\left[ \begin {array} {c} #1 \end {array} \right]}
\newcommand {\onebytwo} [2] {\left[ \begin {array} {c c} #1,& \;#2 \end {array} \right]}
\newcommand {\twobyone} [2] {\left[ \begin {array} {c} #1\\ #2 \end {array} \right]}
\newcommand {\twobytwo} [4] {\left[ \begin {array} {c c} #1& #2\\ #3& #4 \end {array} \right]}

\newcommand {\twobytworight} [4] {\left[ \begin {array} {r r} #1& #2\\#3& #4 \end {array} \right]}

\title {Unattainability of a Perturbation Bound for Indefinite Linear Least Squares Problems}
 
\author{Joseph F. Grcar\thanks{6059 Castlebrook Drive, Castro Valley, CA 94552 USA (jfgrcar@comcast.net).}}

\begin {document}

\maketitle

\begin {abstract} Contrary to an assumption made by Bojanczyk, Higham, and Patel [\textit {SIAM J.\ Matrix Anal.\ Appl.}, 24(4):914--931, 2003], a perturbation bound for indefinite least square problems is capable of arbitrarily large overestimates for all perturbations of some problems. For these problems, the hyperbolic QR factorization algorithm is not proved to be forward stable because the error bound systematically overestimates the solution error of backward stable methods.

\end {abstract}

\begin {keywords}
condition number, forward stability, indefinite linear least squares, perturbation theory
\end {keywords}

\begin {AMS}
primary 65F20; secondary 15A63, 65F35
\end {AMS}


\section {Introduction}

This paper studies a perturbation bound for indefinite linear least squares problems \citep* {Bojanczyk2003}. The problems where the bound is an arbitrarily large overestimate for all perturbations are characterized, and examples of overestimates  are constructed. 

Background material is presented in section \ref {sec:background}. The case of perturbations just to $A$ is considered in section \ref {sec:A}. The case of perturbations to both $A$ and $b$ is considered in section \ref {sec:Ab}.

\section {Background} 
\label {sec:background}

This section describes a canonical form for indefinite least squares problems and explains the perturbation bound.

\subsection {Indefinite Linear Least Squares (ILS)}
 
The ILS problem is
\begin {displaymath}
x = \arg \min_u \, (b - Au)^t J (b - A u)
\end {displaymath}
where $A \in \reals^{m \times n}$, $m > n$, $b \in \reals^m$, $x \in \reals^n$, $J \in \reals^{m \times m}$ is a diagonal matrix with only $\pm 1$ on the diagonal, and $A^t J A$ is positive definite. The positive definite hypothesis reduces the problem to solving the first-order optimality condition, $A^t J (b - A x) = 0$.

Matrix theory for indefinite inner products is developed by \citet* {Gohberg1983}, and for the $J$ inner product by \citet {Higham2003}. Questions in control theory that can be expressed as indefinite linear least squares problems are described by \citet* {Chandrasekaran1998a}. They and \citet {Xu2004} present solution algorithms for the ILS problem that are backward stable.

\subsection {Generalized Singular Value Decomposition}

The following matrix decomposition provides a standard form for ILS problems. Suppose the ILS problem is arranged so that $J$ has positive and negative blocks of orders $m_+$ and $m_-$,
\begin {displaymath}
J = \twobytwo {I} {} {} {-I} \qquad A = \twobyone {A_+} {A_-} \qquad b = \twobyone {b_+} {b_-} \, .
\end {displaymath}
Note $m_+ \ge n$ because $A_+^t A_+ - A_-^t A_-$ is positive definite. The generalized singular value decomposition \citep {Paige1981} then takes the form
\begin {displaymath}
A = {\twobytwo {Q_+} {} {} {Q_-}} \twobyone {C} {S} {X} \qquad \begin {array} {r r c l} \mbox {$m_+ \times n$:}& \; C& =& \diag (\cos (\theta_i)_{\,i=1,\dots,n})\\ \noalign {\smallskip} \mbox {$m_- \times n$:}& S& =& \diag(\sin (\theta_j)_{\,j=1,\dots, \min \{m_-, n\}} ) \end {array}
\end {displaymath}
where $Q_+$, $Q_-$ are square, orthogonal matrices of orders $m_+$, $m_-$, respectively, $X$ is invertible, and the angles satisfy,
\begin {displaymath}
\fourthpi > \theta_1 \ge \theta_2 \ge \cdots \ge \theta_{n} \ge 0 \quad \mbox {and} \quad \mbox {$\theta_i = 0$ for $i > m_-$} \, .
\end {displaymath}
The upper bound on the angles assures $A^t J A = X^t (C^t C - S^t S) X$ is positive definite. 

\subsection {Perturbation Bound}

\newcommand {\maxformula} {\max \left\{  {\| \DA \|_F \over \| \bfA \|_F}, \, {\| \Db \|_2 \over \| \bfb \|_2} \right\}}

Suppose $A + \DA$ and $b + \Db$ are perturbed by $\epsilon$ in a ``normwise relative'' sense,
\begin {equation}
\label {eqn:epsilon}
\left. \begin {array} {r c l} \| \DA \|_F& \le& \epsilon \| \bfA \|_F\\ \noalign {\smallskip} \| \Db \|_2& \le& \epsilon \| \bfb \|_2\end {array} \right\} \quad \mbox {equivalently} \quad
\epsilon = \maxformula ,
\end {equation}
where $\bfA$ and $\bfb$ are normalizing factors, and suppose $(A+\DA)^t J (A+\DA)$ is positive definite. \citet [p.\ 917, eqn.\ 2.8] {Bojanczyk2003} show the change to the solution $x + \Dx$ is at most 
\begin {equation}
\label {eqn:bound}
\begin {array} {r c l}
\displaystyle {\| \Dx \|_2 \over \| x \|_2}& \le& \displaystyle  \| \bfA \|_F  \left( \| (A^t J A)^{-1} A^t \|_2 + \| (A^t J A)^{-1} \|_2 \, {\| r \|_2 \over \| x \|_2}  \right) \epsilon
\\ \noalign {\medskip}
& +& \displaystyle {\| \bfb \|_2 \over \| x \|_2} \, \| (A^t J A)^{-1} A^t \|_2 \, \epsilon + \Oe \, ,
\end {array}
\end {equation}
where $r = b - A x$ is the ILS residual. Bojanczyk et al.\ (p.\ 918, line -10) ``believe that (\ref {eqn:bound}) is nearly attainable.'' 

\section {Perturbations to $\mathbf A$}
\label {sec:A}

For simplicity, this section considers perturbations only to $A$. The first subsection \ref {sec:attainable} uses differential calculus to examine what is meant by an attainable perturbation bound. Subsection \ref {sec:first} presents the simplest ILS problem, for which (\ref {eqn:bound}) can be an overestimate for perturbations to $A$. Perturbations to both $A$ and $b$ are more complicated and are deferred to section \ref {sec:Ab}.

\subsection {Attainable Bounds} 
\label {sec:attainable}

Taylor series are the classic perturbation formulas. For the ILS problem, the function $f(A) = (A^t J A)^{-1} A^t J \, b$ with $b$ fixed gives the dependence of the solution on the matrix, $x = f(A)$. The perturbed value $f(A+\DA)$ is approximated by a power series of order $p$ in $\DA$ under two conditions. (1) $f$ must have continuous partial derivatives of order $p+1$, which is always true because $f$ is algebraic. (2) $f$ must be defined on the line segment from $A$ to $A + \DA$, which amounts to $(A+ s \DA)^t J (A + s \DA)$ being positive definite for all $0 \le s \le 1$. Taylor's formula of order $p=1$ is then \citep [p.\ 245, thm.\ 20.16] {Bartle1964} 
\begin {equation}
\label {eqn:taylor}
f(A + \DA) = f(A) + \underbrace {\Df 1 (A) \, \DA}_{\hspace*{-3em}\displaystyle \JxA \Vector (\DA)\hspace*{-3em}} + \R (\DA, A) \, ,
\end {equation}
where $\Df 1$ is the first derivative of $f$, and $\R$ is the remainder.

The first derivative is represented by a Jacobian matrix. For specificity, the ``$\Vector$'' construction is needed to order matrix entries into column vectors so partial derivatives with respect to matrix entries can be placed into the Jacobian matrix. Thus, $\Vector (A)$ is the column vector of entries $A_{i,j}$ with $(j,i)$ in lexicographic order, and $\JxAvector$ is the Jacobian matrix of partial derivatives of $f$ evaluated at the given $A$. Notation for Jacobian matrices has a subscript for the name of the function (or the dependent variable), followed by the value for the independent variable in parentheses. The ``$\Vector$'' is omitted for brevity; thus $\JxAvector$ is written $\JxA$.

The remainder is bounded by $\| \DA \|_F^2$ times a coefficient independent of $\DA$ only when $\DA$ is restricted in some way.  Let $\Hessian_k$ be the Hessian matrix for the $k$-th entry of $f(A)$. The $k$-th entry of $\R$ is ${\textstyle {1 \over 2}} (\Vector (\DA))^t \, \Hessian_k \Vector (\DA)$ where $\Hessian_k$ is evaluated at $A + s_k \DA$ for some $0 \le s_k \le 1$, which depends on $A$ and $\DA$. Choose $\delta > 0$ so $(A+P)^tJ(A+P)$ is positive definite when $\| P \|_F \le \delta$. Let $\eta_k$ be the maximum of $\| \Hessian_k \|_2$ over  these $A+P$, and let $\eta = (\eta_1^2 + \eta_2^2 + \dots + \eta_n^2)^{1/2}$. Thus $\| \R (P, A) \|_2 \le \halfeta \| P \|_F^2$.

Subtracting $f(A)$ from both sides of equation (\ref {eqn:taylor}) and applying norms gives the general perturbation bound,
\begin {equation}
\label {eqn:differential-bound-for-A-1}
\quad
\| \Dx \|_2 \le \| \JxA \|_2 \, \| \DA \|_F + \halfeta \| \DA \|_F^2  \qquad \mbox {for $\| \DA \|_F \le \delta$} \, ,
\end {equation}
where $\delta$, $\eta$ only depended on $A$, $f$. This bound is remarkable in three respects. 

\textit {First, (\ref {eqn:differential-bound-for-A-1}) is approximately attainable in the following sense.} Let the induced norm $\| \JxA \|_2$ be attained at $\Vector (U)$, which may be scaled so that $\halfeta \| U \|_F \le \| \JxA \|_2$ and $\| U \|_F \le \delta$. For $0 \le t \le 1$, the perturbation $tU$ changes the solution by $\Dx(t) = f(A + t U) - f(A)$, which by (\ref {eqn:taylor}) equals $\JxA \Vector (t U) + \R (tU, A)$. Triangle inequalities then imply 
\begin {equation}
\label {eqn:triangle}
\| \Dx(t) \|_2 \quad \mbox {lies between} \quad \| \JxA \|_2 \, \| t U \|_F \pm \halfeta \| t U \|_F^2 \, .
\end {equation}
Thus, $\DA = t U$ can be chosen so $\| \Dx(t) \|_2$ is arbitrarily close to the bound (\ref {eqn:differential-bound-for-A-1}).

\textit {Second, the coefficient $\| \JxA \|_2$ is the smallest that can appear in bounds of the form (\ref {eqn:differential-bound-for-A-1}).} Suppose there is a comparable bound,
\begin {displaymath}
\| \Dx \|_2 \le \alpha \| \DA \|_F + {\textstyle {\beta \over 2}} \| \DA \|_F^2 \qquad \mbox {for $\| \DA \|_F \le \delta$} \, .
\end {displaymath}
If $\DA = t U$, then the bounds are quadratic polynomials of $t$. This upper bound exceeds the lower bound (\ref {eqn:triangle}) for all $0 \le t \le 1$ if and only if $\alpha \ge \| \JxA \|_2$.

\textit {Third, $\| \JxA \|_2$ properly scaled gives the condition number of $x$ with respect to $A$.} Bound (\ref {eqn:differential-bound-for-A-1}) can be rewritten as
\begin {equation}
\label {eqn:differential-bound-for-A-2}
{\| \Dx \|_2 \over \| x \|_2}  \le \underbrace {\left( {\| \bfA \|_F \over \| x \|_2} \, \| \JxA \|_2 \right)}_{\displaystyle \XxA} {\| \DA \|_F \over \| \bfA \|_F } + \underbrace {\halfeta \| \DA \|_F^2 \over \| x \|_2}_{\displaystyle \O (\| \DA \|_F^2)} \, . 
\end {equation}
There may be many ways to define condition numbers, but since the scaled coefficient in (\ref {eqn:differential-bound-for-A-2}) is the smallest possible, the value of this coefficient must be  the condition number of $x$ with respect to $A$ and the scaled norms, $\XxA$. Jacobian matrices were first used to study conditioning by \citet [p.\ 292, thm.\ 4] {Rice1966}.

It is an open question to find a simple formula for the $\XxA$ of general ILS problems. There is a formula for $\JxA$ which leads to a bound for $\XxA$ that appears in (\ref {eqn:bound}). If $F(A, x) = A^t J (b - A x) \equiv 0$ (optimality condition), then by implicit differentiation $\JxA = - (\J_F (x))^{-1} \J_F (A)$.  In the present case, by inspection $\J_F(x) = - A^t J A$, and  by the product rule for differentiation
\begin {displaymath}
\J_F (A)
=
\setlength {\arraycolsep} {0.25em}
\left[ 
\begin {array} {c c c}
(Jr)^t\\
& \ddots\\
&& \hspace {0.5em} (Jr)^t
\end {array}
\right]
- 
\left[ 
\setlength {\arraycolsep} {0.33em}
\begin {array} {c c c c c}
x_1 A^t J& \cdots& x_n A^t J
\end {array}
\right] 
= I \otimes (Jr)^t - x^t \otimes A^t J \, ,
\end {displaymath}
where $r = b - Ax$, and $\otimes$ is the Kronecker product ($H \otimes K$ has blocks consisting of $K$ scaled by the corresponding entries of $H$). Combining these matrices by the implicit differentiation formula gives\footnote {Equation (\ref {eqn:JxA}) was derived by \citet [p.\ 918, line 13] {Bojanczyk2003} using perturbation algebra.}
\begin {equation}
\label {eqn:JxA}
\JxA = \underbrace {(A^t JA)^{-1} \otimes (Jr)^t}_{\displaystyle \M_1} - \underbrace {x^t \otimes (A^t J A)^{-1} A^t J}_{\displaystyle \M_2} \, .
\end {equation} 
The norms of $\M_1$ and $\M_2$ have simple formulas, which make a bound for $\XxA$,
\begin {equation}
\label {eqn:JxA_bound}
\begin {array} {c}
\| \M_1 \|_2 = \| (A^t J A)^{-1} \|_2 \, \| r \|_2
\qquad
\| \M_2 \|_2 =  \| x \|_2 \, \| (A^t J A)^{-1} A^t \|_2
\\ \noalign {\medskip}
\displaystyle
\XxA = {\| \bfA \|_F \over \| x \|_2} \, \| \M_1 - \M_2 \|_2 \le {\| \bfA \|_F \over \| x \|_2} \, \left( \| \M_1 \|_2 + \|  \M_2 \|_2 \right) = \bhp \, .
\end {array}
\end {equation}
The first line of Bojanczyk et al.'s bound (\ref {eqn:bound}) is $\bhp \kern0.05em \epsilon$.

\subsection {First Example} 
\label {sec:first}

With the choice $\| \bfb \|_2 = 0$, (\ref {eqn:epsilon}) becomes $\Db = 0$ and $\epsilon = \| \DA \|_F / \| \bfA \|_F$, and the first line of (\ref {eqn:bound}) is
\begin {equation}
\label {eqn:bound-for-A}
{\| \Dx \|_2 \over \| x \|_2} \le  \bhp \,  {\| \DA \|_F \over \| \bfA \|_F } + \Oe \, . 
\end {equation}
Comparing (\ref {eqn:differential-bound-for-A-2}) and (\ref {eqn:bound-for-A}), whether (\ref {eqn:bound}) is nearly attainable entails the question of whether the formula $\bhp$ of (\ref {eqn:JxA_bound}) is a good estimate for the condition number $\XxA$.

The answer is ``no'' for the problem with $m = 2$ and $n = 1$,
\begin {displaymath}
A = \twobyone {a_{1,1}} {a_{2,1}} \qquad b = \twobyone 1 1 \qquad J = \twobytworight 1 {} {} {-1} \, .
\end {displaymath}
The solution is
\begin {displaymath}
x = \onebyone {x_1} = (A^t J A)^{-1} A^t J b = \onebyone {\displaystyle {1 \over a_{1,1} + a_{2,1}}} \, ,
\end {displaymath}
and the Jacobian matrix computed by differentiating $x$ is
\begin {displaymath}
\JxA = \onebytwo {\displaystyle {\partial x_1 \over \partial a_{1,1}}} {\displaystyle {\partial x_1 \over \partial a_{2,1}}} = \onebytwo {\displaystyle {-1 \over (a_{1,1} + a_{2,1})^2 }} {\displaystyle {-1 \over (a_{1,1} + a_{2,1})^2}} .
\end {displaymath}
Consider matrices with the following generalized singular value decomposition,
\begin {displaymath}
\left[ \begin {array} {c} a_{1,1}\\ a_{2,1} \end {array} \right] = \left[ \begin {array} {c} \cos (\theta)\\ \sin (\theta) \end {array} \right] \, .
\end {displaymath}
Choose $\| \bfA \|_F = \| A \|_F$ for the remaining scale factor. If $\theta = \fourthpi - \param$ and $0 < \alpha \le 0.3$, then from (\ref {eqn:differential-bound-for-A-2}) $\XxA = \sec (\param) \approx 1.0$, while from (\ref {eqn:bound-for-A}) $\bhp = 2 \csc (2 \param) \approx 1.0/\param$, in which $1.0$ indicates values rounded to two decimal digits for these $\alpha$. The formulas have been obtained with Mathematica \citep {Wolfram2003}, or they are easily calculated by hand. Thus, $\bhp \approx 1/\param$ is not a good approximation to $\XxA \approx 1$. For example, if $\param = 0.01$, then $\XxA = 1.00005$ and $\bhp = 100.007$ rounded to six digits. This example has considerable cancellation in the formula (\ref {eqn:JxA}) for $\JxA$, 
\begin {displaymath}
\begin {array} {r c l}
\| \JxA \|_2 = \| \M_1 - \M_2 \|_2 = 2^{-1/2} \sec^2 (\alpha)& \approx& 2^{-1/2} \, ,
\\ \noalign {\smallskip}
\| \M_1 \|_2 = \| \M_2 \|_2 = 2^{-1/2} \csc (2 \param) \sec (\param)& \approx& 2^{-3/2} / \param \, .
\end {array}
\end {displaymath}

\section {Perturbations to $\mathbf A$ and $\mathbf b$}
\label {sec:Ab}

This section shows that (\ref {eqn:bound}) is not nearly attainable in general for perturbations to both $A$ and $b$. The first subsection \ref {sec:analysis} examines the types of bounds that are attainable for problems that depend on a matrix and a vector. Subsection \ref {sec:theorem} characterizes the situations in which (\ref {eqn:bound}) is an overestimate, then subsection \ref {sec:second} presents an example.

\subsection {Attainable Bounds Redux}
\label {sec:analysis}

An attainable perturbation bound for chan\-ges only to $b$ can be constructed as in section \ref {sec:attainable}. The partial derivatives of $x = (A^t J A)^{-1} A^t J b$ with respect to $b$ are  the entries of the coefficient matrix $(A^t J A)^{-1} A^t J$, which is then $\Jxb$. Thus, like equation (\ref {eqn:differential-bound-for-A-2}) for $A$, so for $b$,
\begin {equation}
\label {eqn:differential-bound-for-b}
{\| \Dx \|_2 \over \| x \|_2}  \le \underbrace {\left( {\| \bfb \|_2 \over \| x \|_2} \, \| \Jxb \|_2 \right)}_{\displaystyle \Xxb} {\| \Db \|_2 \over \| \bfb \|_2 } \quad \mbox {where also} \quad \Xxb = {\| \bfb \|_2 \, \| \M_2 \|_2 \over \| x \|_2^2} \, . 
\end {equation} 
Note for future reference that $\Xxb$ has a formula in terms of the matrix $\M_2$ in the formulas (\ref {eqn:JxA_bound}) for $\XxA$ and $\bhp$. There is no higher-order remainder term because $x$ is a linear function of $b$. As a result, this bound can be exactly attained. 

When both $A$ and $b$ can be perturbed, then there are several perturbation bounds. The solution $x = g(A,b) = (A^t J A)^{-1} A^t J b$ has a Taylor series analogous to (\ref {eqn:taylor}),
\begin {equation}
\label {eqn:taylor2}
g(A + \DA, b+\Db) = g(A,b) + \underbrace {\Dg 1 (A,b) \, (\DA, \Db)}_{\hspace*{-3em}\displaystyle \JxAb \Vector (\DA, \Db)\hspace*{-3em}} + \, \calS \, , 
\end {equation}
where $\Dg 1$ is the first derivative of $g$, and $\calS$ is the remainder. $\calS$ is again a ``second order'' expression of $\DA$, $\Db$, but for brevity a bound is not developed here. The ``$\Vector$'' construction is extended to $\Vector (A,b)$ by placing the entries of $b$ after $A$. This ordering makes $\JxAb$ into a $2 \times 1$ block matrix of the separate Jacobian matrices with respect to $A$ and $b$, $\JxAb = \onebytwo {\JxA} {\Jxb}$.

The simplest way to convert (\ref {eqn:taylor2}) into a bound applies the $2$-norm to both sides and uses the triangle inequality on the right to treat each block of $\JxAb$ separately,
\begin {eqnarray}
\label {eqn:separate_bound}
\qquad {\| \Dx \|_2 \over \| x \|_2} &\le& \underbrace {\| \JxA \|_2 \, \| \bfA \|_F \over \strut \| x \|_2}_{\displaystyle \XxA} \, {\| \DA \|_F \over \| \bfA \|_F} + \underbrace {\|  \Jxb \|_2 \, \| \bfb \|_2 \over \strut \| x \|_2}_{\displaystyle \Xxb} \, {\| \Db \|_2 \over \| \bfb \|_2} + \Sremainder
\\ \noalign {\smallskip}
\label {eqn:sum_bound}
& \le& \left[ \XxA + \Xxb \right] \epsilon + {\| \calS \|_2 / \| x \|_2} \, .
\end {eqnarray}
Like (\ref {eqn:differential-bound-for-A-2}, \ref {eqn:differential-bound-for-b}), bound (\ref {eqn:separate_bound}) also can be approximately attained.

A more complicated way to convert (\ref {eqn:taylor2}) into a bound requires a norm for $\Vector (\DA, \Db)$. Such a norm is $\| \Vector (\DA, \Db) \|_{\rm max} = \epsilon$ of (\ref {eqn:epsilon}). The compatible norm for $2 \times 1$ block matrices has an upper bound,
\begin {displaymath}
\begin {array} {l}
\displaystyle \makebox [7em] {\hfill $\| \JxAb \|$} = \max_{\DM,\Dv} {\| \JxAb \Vector (\DM, \Dv) \|_2 \over \| \Vector (\DM, \Dv) \|_{\rm max}}
\\ \noalign {\medskip}
\displaystyle \makebox [5em] {} \le \max_{\DM,\Dv} {\| \JxA \|_2 \, \| \bfA \|_F {\| \DM \|_F \over \| \bfA \|_F} + \| \Jxb \|_2 \, \| \bfb \|_2 {\| \Dv \|_2 \over \| \bfb \|_2} \over \| \Vector (\DM, \Dv) \|_{\rm max}}
\\ \noalign {\medskip}
\displaystyle \makebox [4em] {} \le \; \| \JxA \|_2 \, \| \bfA \|_F +  \| \Jxb \|_2 \, \| \bfb \|_2 \, .
\end {array}
\end {displaymath}
Applying norms to (\ref {eqn:taylor2}) and using this matrix norm and its bound gives two bounds,\footnote {\citet [p.\ 918, line -15] {Bojanczyk2003} use the symbol $\psi = \XxA + \Xxb$.}
\begin {equation}
\label {eqn:joint_bound}
\qquad {\| \Dx \|_2 \over \| x \|_2} \le \underbrace {\| \JxAb \| \over \| x \|_2}_{\displaystyle \XxAb} \epsilon + \Sremainder \le \underbrace {\| \JxA \|_2 \, \| \bfA \|_F +  \| \Jxb \|_2 \, \| \bfb \|_2 \over \| x \|_2}_{\displaystyle \XxA + \Xxb} \epsilon + \Sremainder \, . 
\end {equation}
The first of these bounds can be approximately attained in the sense of (\ref {eqn:triangle}).

Altogether (\ref {eqn:separate_bound}, \ref {eqn:joint_bound}) and (\ref {eqn:sum_bound}) with (\ref {eqn:JxA_bound}) give four bounds with $\epsilon$ as in (\ref {eqn:epsilon}),
\begin {equation}
\label {eqn:all_bounds}
\qquad
\framebox {
\begin {minipage} {21.25em}
$\displaystyle
{\| \Dx \|_2 \over \| x \|_2} \le 
\left\{ \begin {array} {r}
\displaystyle \plethorasup {i} \; \XxA {\| \DA \|_F \over \| \bfA \|_F} + \Xxb {\| \Db \|_2 \over \| \bfb \|_2} 
\\ \noalign {\medskip}
\plethorasup {ii} \; \XxAb \, \epsilon \,
\end {array} \right\} 
$
\end {minipage}%
\hspace* {-0.5em}%
\raisebox {-4.2ex} {
\begin {minipage} [b] {11em}
$\displaystyle \begin {array} {c}
\le \plethorasup {iii} \left[ \XxA + \Xxb \right] \epsilon
\\ \noalign {\smallskip}
\le \plethorasup {iv} \underbrace {\left[ \bhp+ \Xxb \right] \epsilon}_{\mbox {\hspace{-0.25em}bound (\ref {eqn:bound})\hspace{-0.25em}}}  \, ,
\end {array}$
\end {minipage}
}}
\end {equation}
in which the $\O (\epsilon^2)$ remainder term $\| S \|_2 / \| x \|_2$ is omitted for clarity. The two smallest \plethora {i, ii} can be approximately attained, but they are incommensurate. If the separate condition numbers $\XxA$, $\Xxb$ are of different sizes, then the scaled perturbations $\DA$, $\Db$ can be chosen to make the sum of products \plethora {i} much smaller than the product of the sum and maximum \plethora {iii}. However, the bound \plethora {ii} with one condition number is always close to the bound \plethora {iii}, because\footnote {These inequalities are remarked by \citet [p.\ 918, line -13] {Bojanczyk2003} and are proved by \citet [p.\ 2937, eqn.\ 2.10] {Grcar2010f}.} 
\begin {displaymath}
1 \le {\XxA + \Xxb \over \XxAb} \le 2 \, .
\end {displaymath}
Thus, for some perturbations, bound \plethora {i} without $\epsilon$ may be far smaller than all the other bounds with $\epsilon$.

\subsection {When (\ref {eqn:bound}) is not Attainable}
\label {sec:theorem}

For perturbations only to $A$, the discussion of section \ref {sec:attainable})  shows a bound is not attainable (to within high order terms) only if the coefficient in the bound overestimates the condition number, and (\ref {eqn:JxA_bound}) shows $\bhp \ge \XxA$ is caused by cancellation in $\| \M_1 - \M_2 \|_2$ compared to $\| \M_1 \|_2 + \| \M_2 \|_2$. 

For perturbations to both $A$ and $b$, (\ref {eqn:all_bounds}) shows (\ref {eqn:bound}) is not attainable but it is close to an attained bound unless $\bhp + \Xxb$ overestimates $\XxA + \Xxb$. This situtation is much more complicated because $\Xxb$ my be comparable to or larger than the other terms. Theorem {\ref {thm:conditions} is a precise statement of when an overestimate occurs. Some preparation is needed to prove the theorem.

\begin {lemma}
\label {lem:m}
If $m_1$ and $m_2$ are positive, then
\begin {displaymath}
{2 \over 1 + \mu} \le^{(a)} {m_1 + m_2 \over m_2} \le^{(b)} {2 \over 1 - \mu} \quad \mbox {where} \quad \mu = {| m_1 - m_2 | \over m_1 + m_2} \, .
\end {displaymath}
\end {lemma}
\begin {proof} 
Taking reciprocals, multiplying by $2(m_1+m_2)$, and subtracting $m_1 + m_2$, shows the inequalities are equivalent to $| m_1 - m_2 | \ge m_2 - m_1 \ge - | m_1 - m_2 |$.
\end {proof}

\begin {lemma}
\label {lem:uvw}
If $u$, $v$, $w$ are positive and 
\begin {displaymath}
\rho = {w + u \over v + u} \ge 1 \, , \quad \mbox {then} \quad {w \over u} + 1 \ge^{(a)} \rho \ge^{(b)} {1\over 2} \left( {w \over u} + 1\right) \quad \mbox {and} \quad {w \over v} >^{(c)} \rho >^{(d)} {1 \over 2} {w \over v} \, ,
\end {displaymath}
where at least one of (b) and (d) is true, but not necessarily both.
\end {lemma}
\begin{proof}
(a) $\rho \le (w+u)/u=(w/u)+1$.
(c) $w + u = \rho v + \rho u$ so $w - \rho v = (\rho-1) u \ge 0$, hence $w \ge \rho v$, and then $w/v \ge \rho$. This is the only use of $\rho \ge 1$.
(b) If $u \ge v$, then $\rho \ge (w + u) / (2 u) = ((w/u) + 1)/2$. 
(d) If $v \ge u$, then $\rho \ge (w + u) / (2v) \ge w/(2v)$.
\end{proof}

\begin {theorem}
\label {thm:conditions}
If
\begin {displaymath}
\rho = {\bhp + \Xxb \over \XxA + \Xxb} 
\qquad
\lambda_1 = {\| \M_1 \|_2 + \| \M_2 \|_2 \over \| \M_1 - \M_2 \|_2}
\qquad
\lambda_2 = {\| \bfA \|_F  \| x \|_2 \over \| \bfb \|_2}
\end {displaymath}
where $\bhp$, $\M_1$, $\M_2$, $\XxA$, $\Xxb$ are given by (\ref {eqn:JxA_bound}, \ref {eqn:differential-bound-for-b}) and no terms vanish, then
\begin {displaymath}
\framebox {$\strut
\hspace*{0.5em}
2 \rho \ge^* \lambda_1 \ge \rho
\quad \mbox {and} \quad
\rho + {1 \over 2} \ge^* \lambda_2 \ge  {\rho \over 2} - 1
\hspace*{0.5em}
$}
\end {displaymath}
in which at least one of the inequalities (*) is true, but not necessarily both. Consequently, bound (\ref {eqn:bound}) is a uniformly large overestimate for all $\DA$ and $\Db$ if and only if both $\lambda_1$ and $\lambda_2$ are much larger than $1$, in which case the factor of overestimation is at least either $\lambda_1 2^{-1}$ or $\lambda_2 - 2^{-1}$.
\end {theorem}

\begin {proof}
From (\ref {eqn:JxA_bound}) note $\rho \ge 1$. With (\ref {eqn:differential-bound-for-b}) and some rearrangement,
\begin {displaymath}
\rho = {(\| \M_1 \|_2 + \| \M_2 \|_2) + (\| \M_2 \|_2 \, \| \bfb \|_2) / ( \| \bfA \|_F \, \| x \|_2) \over (\| \M_1 - \M_2 \|_2) + (\| \M_2 \|_2 \, \| \bfb \|_2) / ( \| \bfA \|_F \, \| x \|_2)} = {w + u \over v + u} \, ,
\end {displaymath}
where
\begin {displaymath}
u = (\| \M_2 \|_2 \, \| \bfb \|_2) / ( \| \bfA \|_F \, \| x \|_2) \quad
v = \| \M_1 - \M_2 \|_2 \quad
w = \| \M_1 \|_2 + \| \M_2 \|_2 \, .
\end {displaymath}

Part 1, inequality $\lambda_1 \ge \rho$: Note $\lambda_1 = w/v$, and from lemma \ref {lem:uvw} (c) $w/v \ge \rho$.

Part 2, inequality $\lambda_2 \ge {\rho \over 2} - 1$: 
Let $m_i = \| \M_i \|_2$, and let $\mu = \left| m_1 - m_2 \right| / (m_1 + m_2)$. 
From one of the triangle inequalities, $1 / \lambda_1 \ge \mu$; therefore from part 1 of this proof, (2.i) $1/ \rho \ge \mu$.
From lemma \ref {lem:m} (b), $2/(1-\mu) \ge (m_1 + m_2)/m_2$, which combines with (2.i) to give, (2.ii) $2/(1-\rho^{-1}) \ge (m_1 + m_2)/m_2$.
Multiplying (2.ii) by $\lambda_2$ and noting $\lambda_2 (m_1+m_2)/m_2 = w/u$ implies, (2.iii) $2\lambda_2/(1-\rho^{-1}) \ge w/u$. Lemma \ref {lem:uvw} (c) is (2.iv) $w/u \ge \rho - 1$. 
The last two inequalities (2.iii, 2.iv) combine to $2\lambda_2/(1-\rho^{-1}) \ge \rho - 1$, which weakens to $\lambda_2 \ge {\rho \over 2} - 1$.

Part 3, inequality $2 \rho \ge^* \lambda_1$: 
If lemma \ref {lem:uvw} (d) is true, then $\rho \ge (w/v)/2$ where as noted, $w/v = \lambda_1$. Rearranging gives $2 \rho \ge \lambda_1$.

Part 4, inequality $\rho + {1 \over 2} \ge^* \lambda_2$:
If lemma \ref {lem:uvw} (b) is true, then $2 \rho -1 \ge w/u$ where again $w/u= \lambda_2 (m_1+m_2)/m_2$ so altogether, (4.i) $2 \rho -1 \ge  \lambda_2 (m_1+m_2)/m_2$.
From lemma \ref {lem:m} (a) $(m_1 + m_2)/m_2 \ge 2/(1+\mu)$, which combines with (2.i) to give (4.ii) $(m_1 + m_2)/m_2 \ge 2/(1+\rho^{-1})$.
Multiplying (4.ii) by $\lambda_2$ and combining with (4.i) gives $2\rho-1 \ge 2\lambda_2/(1+\rho^{-1})$.
The latter rearranges to $\rho + (1 - \rho^{-1}) / 2 \ge \lambda_2$ which weakens by discarding the negative term.
\end {proof}

The formulas of theorem \ref {thm:conditions} resemble some given by \citet [p.\  917] {Bojanczyk2003}. They introduce three quantities 
$E_1 = \| \M_2 \|_2  \, \| \bfb \|_2 / \| x \|_2$, $E_2 = \| \M_2 \|_2 \, \| \bfA \|_F$, and $E_3 = \| \M_1 \|_2 \, \| \bfA \|_F$, and they comment: ``the bound (\ref {eqn:bound}) can fail to be achieved for some $\Db$ and $\DA$ only if $E_1 < E_2 \approx E_3$ and there is substantial cancellation in the expression $-(A^tJA)^{-1} A^t J \DA x + (A^tJA)^{-1} \DA^t Jr$ for all $\DA$. We can show in various special cases that these circumstances cannot arise, but we have been unable to establish attainability of the bound (\ref {eqn:bound}) in general.'' 

\subsection {Second Example}
\label {sec:second}

The situation predicted by theorem \ref {thm:conditions} does arise for ILS problems with the following generalized singular value decomposition for $A$,
\begin {displaymath}
A = \left[ \begin {array} {c c c} \cos(\fourthpi-\param)&\;& 0\\ 0&& 1\\ \sin (\fourthpi-\param)
&& 0 \end {array} \right] \twobytwo {1} {0} {0} {\;\param
^{-1}} \qquad b = \left[ \begin {array} {r} 1\\ 1\\ 1 \end {array} \right] \qquad J = \left[ \begin {array}
{c r c} 1\\ & \kern0.6em 1\\ && -1 \end {array} \right] 
\end {displaymath}
where $0 < \param \le \fourthpi$. The same parameter is used in both columns of $A$ for simplicity. Choose the customary scale factors $\bfA = A$ and $\bfb = b$. The following expressions have again been derived with Mathematica ($\| \JxA \|_2$ is beyond easy hand calculation),
\begin {displaymath}
\begin {array} {c}
\displaystyle
\| \bfA \|_F \approx \param^{-1} \qquad
\| \bfb \|_2 \approx 3^{1/2} \qquad
\| r \|_2 \approx 2^{-1/2} \qquad
\| x \|_2  \approx 2^{-1/2}\\ \noalign {\medskip}
\| (A^tJA)^{-1} \|_2 \approx (2 \param)^{-1} \qquad
\| (A^tJA)^{-1} A^t \|_2 \approx (2 \param)^{-1}  \qquad
\| \JxA \|_2 \approx 2^{-1} 3^{1/2}\\ \noalign {\medskip}
\| \M_1 \|_2 \approx 2^{-3/2} \param^{-1} \qquad
\| \M_2 \|_2 \approx 2^{-3/2} \param^{-1} \qquad
\| \M_1 - \M_2 \|_2 \approx 2^{-1} 3^{1/2} \, .
\end {array}
\end {displaymath}
In these formulas, ``$\approx$'' indicates the asymptotic value as $\param$ nears $0$. This example illustrates theorem \ref {thm:conditions},
\begin {displaymath}
\begin {array} {c}
\XxA \approx 2^{-1/2} 3^{1/2} \param^{-1} \qquad
\bhp \approx \param^{-2} \qquad 
\Xxb \approx 2^{-1/2} 3^{1/2} \param^{-1}
\\ \noalign {\medskip}
\rho \approx 6^{-1/2} \param^{-1} \qquad 
\lambda_1 \approx 2^{1/2} 3^{-1/2} \param^{-1} \qquad
\lambda_2 \approx 6^{-1/2} \param^{-1} \, .
\end {array}
\end {displaymath}
Thus, bound (\ref {eqn:bound}) overestimates bound \plethora {iii} by a factor of $(\sqrt 6 \, \param)^{-1}$ asymptotically as $\param$ nears $0$; therefore for all $\DA$ and $\Db$, bound (\ref {eqn:bound}) overestimates $\Dx$ by at least the same factor. For example, if $\param = 0.01$, then $\XxA + \Xxb = 244.937$ and $\bhp + \Xxb = 10123.1$ rounded to six digits. The conditions of theorem \ref {thm:conditions} are,
\begin {displaymath}
\begin {array} {c c c c c}
82.66& \, \ge^* \,& 81.66& \, \ge \,& 41.33\\ \noalign {\smallskip}
2 \rho&& \lambda_1&& \rho
\end {array}
\qquad
\begin {array} {r c c c l}
41.83& \, \ge^* \, & 40.83& \, \ge \,& 19.66\\ \noalign {\smallskip}
\rho + {1 \over 2}&& \lambda_2&& {\rho \over 2} - 1
\end {array}
\end {displaymath}
in which the values have been rounded to two fractional digits, and only one of the inequalities (*) need be true.

\section {Discussion}

The example of section \ref {sec:second} demonstrates that the rounding error analysis of \citet {Bojanczyk2003} does not prove the forward stability of the hyperbolic QR factorization method to solve ILS problems.  Their ``overall conclusion'' (p.\ 929, line -12 and fol.)\ is that the algorithm calculates $x + \Dx$ which satisfies (\ref {eqn:bound}) for $\epsilon = {\mathcal O} (mn {\bf u})$, where $\bf u$ is the unit roundoff.  ``If we make the reasonable assumption that the perturbation bound (\ref {eqn:bound}) is approximately attainable, then our rounding error analysis has shown that the \dots\ [method] is forward stable.'' Forward stable means the algorithm satisfies the same normwise perturbation bounds as backward stable methods \citep [p.\ 130] {Higham2002}. The complicated example of section \ref {sec:second} shows (\ref {eqn:bound}) overestimates $\| \Dx \|$ by an arbitrarily large factor for all perturbations $\DA$, $\Db$ of some ILS problems, in which case backward stable algorithms do satisfy a bound with a much smaller coefficient $\XxA + \Xxb \ll \bhp + \Xxb$. 

Nevertheless, it may be said that the hyperbolic QR factorization method is \textit {provisonally forward stable\/}, in the following sense. If one or both of the quantities $\lambda_1$, $\lambda_2$ of theorem \ref {thm:conditions} are close to $1$, then (\ref {eqn:bound}) approximates an attainable bound, and the method is forward stable for that problem. At least the size of $\lambda_2$ is easy to check: note, $\lambda_2 \approx 1$ is sufficient, but not necessary, for (\ref {eqn:bound}) to be nearly attainable.

\section {Conclusion}

Of greater significance than the stability of any one algorithm is that the simple example of section \ref {sec:first}  shows the ILS problem is among those for which typical normwise formulas (resembling others in numerical linear algebra) are inadequate to estimate the conditioning of the problem. \citet {Malyshev2003} has identified many problems related to least squares that have plausible error bounds that can be far from sharp. It is an open question how to find simple formulas for the condition numbers of these problems. 

\section* {Acknowledgements}

I thank the editor, the reviewers, and Profs.\ Bojanczyk and Higham for corrections, remarks, and suggestions that  improved this paper.


{
\frenchspacing
\raggedright

}

\end {document}